\documentclass[12pt]{article}
\usepackage{amsmath,amsthm,amssymb}

\title{Asymptotical Bounds for Complete Elliptic Integrals of the Second Kind}
\author{\small Miao-Kun Wang, Yu-Ming Chu}
\date{}

\begin{document}

\maketitle
\renewcommand{\thefootnote}{\fnsymbol{footnote}}

{\footnotetext{\hspace*{-5mm}
\begin{tabular}{@{}r@{}p{13cm}@{}}
&M.-K. Wang\\
&Department of Mathematics, Huzhou Teachers College, Huzhou 313000,
China\\
&e-mail: wmk000@126.com\\
\\
&Y.-M. Chu (Corresponding author)\\
&Department of Mathematics, Huzhou Teachers College, Huzhou 313000,
China\\
&e-mail: chuyuming@hutc.zj.cn; tel: +86 572 2321510; fax: +86 572
2321165
\end{tabular}}}

\medskip
\noindent{\bf Abstract}:  In this paper, we establish several
asymptotical bounds for the complete elliptic integrals of the
second kind $\mathcal{E}(r)$, and improve the well-known conjecture
$\mathcal{E}(r)>\pi[(1+(1-r^2)^{3/4})/2]^{2/3}/2$ for all
$r\in(0,1)$ proposed by M. Vuorinen.



\bigskip
\bigskip
\noindent{\bf\large 1. Introduction}
\bigskip
\setcounter{section}{1}
\setcounter{equation}{0}

For $0<r<1$, the well-known complete elliptic integrals of the first
and second kinds [15, 16] are defined by
\begin{equation*}
{\mathcal{K}}={\mathcal{K}}(r)=\int_{0}^{\pi/2}(1-r^2\sin^{2}\theta)^{-1/2}d\theta,\quad
\mathcal{K}(0^+)=\pi/2,\quad \mathcal{K}(1^{-})=+\infty
\end{equation*}
and
\begin{equation*}
\mathcal{E}=\mathcal{E}(r)=\int_{0}^{\pi/2}(1-r^2\sin^{2}\theta)^{1/2}d\theta,\quad
\mathcal{E}(0^+)=\pi/2,\quad \mathcal{E}(1^{-})=1,
\end{equation*}
respectively.

It is well known that the complete elliptic integrals have many
important applications in physics, engineering, geometric function
theory, quasiconformal analysis, theory of mean values, number
theory and other related fields [4-6, 9, 16, 24, 26, 27].

Recently, the complete elliptic integrals have attracted the
attention of numerous mathematicians. In particular, many remarkable
properties and inequalities for the complete elliptic integrals can
be found in the literature [1-4, 7, 8, 10, 12-14, 19, 21, 22, 26,
29-33].

Let $L(1,r)$ be the arc length of an ellipse with semiaxes 1 and
$r\in(0,1)$, then
\begin{equation*}
L(1,r)=4\mathcal{E}(r'),
\end{equation*}
where $r'=\sqrt{1-r^2}$, and in what follows we use the symbol $r'$
to stand for $\sqrt{1-r^2}$ for $r\in[0,1]$. Inspired by a result of
T. Muri [23], who pointed out in 1883 that $L(1,r)$ can be
approximated by $2\pi[(1+r^{3/2})/2]^{2/3}$, M. Vuorinen [28]
conjectured that inequality
\begin{equation}
{\mathcal{E}}(r)>\frac{\pi}{2}\left(\frac{1+{r'}^{3/2}}{2}\right)^{2/3}
\end{equation}
holds for all $r\in(0,1)$. The conjecture was proved by R. W.
Barnared et al. in [13]. Later, they also provided an upper bound
for $\mathcal{E}(r)$ [14]
\begin{equation}
{\mathcal{E}}(r)<\frac{\pi}{2}\left(\frac{1+{r'}^2}{2}\right)^{1/2},\quad
0<r<1.
\end{equation}

In [2], H. Alzer and S.-L. Qiu announced that (without proof) the
inequality
\begin{equation}
\mathcal{E}(r)<\frac{\pi}{4}\left(\sqrt{1-\overline{\alpha
}r^2}+\sqrt{1-\overline{\beta}r^2}\right)
\end{equation}
holds for all $r\in(0,1)$ with $\overline{\alpha}=1/2-\sqrt{2}/4$
and $\overline{\beta}=1/2+\sqrt{2}/4$, and stated that inequality
(1.3) was offered by one of the referees. Later, H. Kazi and E.
Neuman [22] gave a full proof of inequality (1.3) based on the
two-point Gauss-Chebyshev quadrature formula with the remainder
[11].

The purpose of this paper is to establish several asymptotical upper
and lower bounds for complete elliptic integrals of the second kind,
which improve inequalities (1.1) and (1.2), and give an lower bound
of inequality (1.3). Our main results are the following Theorems 1.1
and 1.2.

\medskip
{\bf Theorem 1.1.} If $\alpha,\beta\in(0,1/2]$, then the double
inequality
\begin{align}
\frac{\pi}{4}\left(\sqrt{\beta+(1-\beta){r'}^2}+\sqrt{(1-\beta)+\beta{r'}^2}\right)<\mathcal{E}(r)\nonumber\\
<\frac{\pi}{4}\left(\sqrt{\alpha+(1-\alpha){r'}^2}+\sqrt{(1-\alpha)+\alpha{r'}^2}\right)
\end{align}
holds for all $r\in(0,1)$ if and only if $\beta\leq
1/2-2\sqrt{2({\pi}^2-8)}/{\pi}^2$ and $\alpha\geq 1/2-\sqrt{2}/4$.

\medskip
{\bf Theorem 1.2.} If $t_{1},t_{2}\in[1/2,1]$ and $p\in[1/2,2]$,
then the double inequality
\begin{align}
&2^{p-2}\pi(1+r')^{1-2p}\left\{\left[t_{1}+(1-t_{1})r'\right]^2+\left[(1-t_{1})+t_{1}r'\right]^2\right\}^{p}<\mathcal{E}(r)\nonumber\\
&<2^{p-2}\pi(1+r')^{1-2p}\left\{\left[t_{2}+(1-t_{2})r'\right]^2+\left[(1-t_{2})+t_{2}r'\right]^2\right\}^{p}
\end{align}
holds for all $r\in(0,1)$ if and only if $t_{1}\leq
1/2+\sqrt{1/(4p)}/2 $ and $t_{2}\geq 1/2+\sqrt{(4/\pi)^{1/p}-1}/2$.

\bigskip
\bigskip
\noindent{\bf\large 2. Lemmas}
\bigskip
\setcounter{section}{2} \setcounter{equation}{0}

In order to prove our main results we need several formulas and
lemmas, which we present in this section.
\medskip

For $r\in(0,1)$, the following formulas were presented in [9,
Appendix E, pp. 474-475, and Theorem 3.21(7)].
\begin{equation*}
\frac{d\mathcal{K}}{dr}=\frac{\mathcal{E}-{r'}^2\mathcal{K}}{r{r'}^2},\quad
\frac{d\mathcal{E}}{dr}=\frac{\mathcal{E}-\mathcal{K}}{r},\quad
\frac{d({\mathcal{E}}-{r'}^2{\mathcal{K}})}{dr}=r{\mathcal{K}},
\end{equation*}
\begin{equation*}
\frac{d({\mathcal{K}}-{\mathcal{E}})}{dr}=\frac{r{\mathcal{E}}}{{r'}^2},
\quad
\mathcal{E}\left(\frac{2\sqrt{r}}{1+r}\right)=\frac{2\mathcal{E}-{r'}^2\mathcal{K}}{1+r},
\quad \lim_{r\rightarrow 1} {r'}^{\alpha}\mathcal{K}(r)=0 \
(\alpha>0).
\end{equation*}

\medskip
{\bf Lemma 2.1} (See [9, Theorem 1.25]). For $-\infty<a<b<\infty$,
let $f,g:[a,b]\rightarrow{\mathbb{R}}$ be continuous on $[a,b]$, and
be differentiable on $(a,b)$, let $g'(x)\neq 0$ on $(a,b)$. If
$f^{\prime}(x)/g^{\prime}(x)$ is increasing (decreasing) on $(a,b)$,
then so are
$$\frac{f(x)-f(a)}{g(x)-g(a)}\ \ \mbox{and}\ \ \frac{f(x)-f(b)}{g(x)-g(b)}.$$
If $f^{\prime}(x)/g^{\prime}(x)$ is strictly monotone, then the
monotonicity in the conclusion is also strict.

\medskip
{\bf Lemma 2.2.} Let $r\in(0,1)$, then

(1) the function $r\rightarrow(\mathcal{E}-{r'}^2\mathcal{K})/r^2$
is strictly increasing from $(0,1)$ onto $(\pi/4,1)$;

(2) the function $r\rightarrow \mathcal{E}/{{r'}^{1/2}}$ is strictly
increasing from $(0,1)$ onto $(\pi/2,\infty)$;

(3) the function $r\rightarrow
(\mathcal{K}-\mathcal{E})/(r^2\mathcal{K})$ is strictly increasing
from $(0,1)$ onto $(1/2,1)$;

(4) the function
$r\rightarrow(\mathcal{E}-{r'}^2\mathcal{K})/(r^2\mathcal{K})$ is
strictly decreasing from $(0,1)$ onto $(0,1/2)$;

(5) the function $r\rightarrow
{r'}^{3/4}(\mathcal{K}-\mathcal{E})/r^2$ is strictly decreasing from
$(0,1)$ onto $(0,\pi/4)$;

(6) the function $r\rightarrow
(\mathcal{E}-{r'}^2\mathcal{K})^2/({\mathcal{E}}^2-{r'}^2{\mathcal{K}}^2)$
is strictly decreasing from $(0,1)$ onto $(1,2)$;

(7) the function
$r\rightarrow[{4(2{\mathcal{E}}-{r'}^2{\mathcal{K}})^2-{\pi}^2}]/r^2$
is strictly increasing from $(0,1)$ onto $({\pi}^2/2,16-{\pi}^2)$.

\medskip
{\bf Proof.}  Parts (1)-(5) can be found in [9, Theorem 3.21(1) and
(8), Exercises 3.43(32) and (46)] and [2, Theorem 15].

For part (6), let $I_{1}(r)=(\mathcal{E}-{r'}^2\mathcal{K})^2$,
$I_{2}(r)={\mathcal{E}}^2-{r'}^2{\mathcal{K}}^2$ and
$I(r)=(\mathcal{E}-{r'}^2\mathcal{K})^2/({\mathcal{E}}^2-{r'}^2{\mathcal{K}}^2)$.
Then $I(r)=I_{1}(r)/I_{2}(r)$, $I_{1}(0)=I_{2}(0)=0$ and
\begin{equation}
\frac{{I_{1}}'(r)}{{I_{2}}'(r)}=\left(\frac{r^2\mathcal{K}}{\mathcal{K}-\mathcal{E}}\right)^2\cdot\frac{\mathcal{E}-{r'}^2\mathcal{K}}{r^2\mathcal{K}}.
\end{equation}

It follows from  parts (3) and (4) together with Lemma 2.1 and (2.1)
that $I(r)$ is strictly decreasing in $(0,1)$. By l'H\^{o}pital's
rule, we have $I(0^+)=2$, while $I(1^-)=1$ is clear.

For part (7), let
$J_{1}(r)={4(2{\mathcal{E}}-{r'}^2{\mathcal{K}})^2-{\pi}^2}$,
$J_{2}(r)=r^2$ and
$J(r)=[{4(2{\mathcal{E}}-{r'}^2{\mathcal{K}})^2-{\pi}^2}]/r^2$. Then
$J(r)=J_{1}(r)/J_{2}(r)$, $J_{1}(0)=J_{2}(0)=0$ and
\begin{equation}
\frac{{J_{1}}'(r)}{{J_{2}}'(r)}=\frac{4(2{\mathcal{E}}-{r'}^2{\mathcal{K}})({\mathcal{E}}-{r'}^2{\mathcal{K}})}{r^2}.
\end{equation}

It follows from part (1) and (2.2) together with Lemma 2.1 that
$J(r)$ is strictly increasing in $(0,1)$. By l'H\^{o}pital's rule,
we get $J(0^+)={\pi}^2/2$, while $J(1^-)=16-{\pi}^2$ is clear.
$\Box$

\medskip
{\bf Lemma 2.3.} Let $r\in(0,1)$ and
\begin{equation*}
g(r)=\frac{(\mathcal{K}-\mathcal{E})(\mathcal{E}-{r'}^2\mathcal{K})+\mathcal{E}[(\mathcal{K}-\mathcal{E})-(\mathcal{E}-{r'}^2\mathcal{K})]}
{(\mathcal{E}-{r'}^2\mathcal{K})^2}.
\end{equation*}
Then $g(r)$ is strictly increasing from $(0,1)$ onto $(3/2,\infty)$.

\medskip
{\bf Proof.} Differentiating $g(r)$ and elaborated computations lead
to
\begin{align}
g'(r)&=\frac{\frac{r\mathcal{E}}{{r'}^2}(\mathcal{E}-{r'}^2\mathcal{K})+(\mathcal{K}-\mathcal{E})r\mathcal{K}+\frac{\mathcal{E}-\mathcal{K}}{r}
[\mathcal{K}-\mathcal{E}-(\mathcal{E}-{r'}^2\mathcal{K})]+\mathcal{E}(\frac{r\mathcal{E}}{{r'}^2}-r\mathcal{K})}
{(\mathcal{E}-{r'}^2\mathcal{K})^2}\nonumber\\
&-\frac{2r\mathcal{K}\left\{(\mathcal{K}-\mathcal{E})(\mathcal{E}-{r'}^2\mathcal{K})+\mathcal{E}[(\mathcal{K}-\mathcal{E})-(\mathcal{E}-{r'}^2\mathcal{K})]\right\}}
{(\mathcal{E}-{r'}^2\mathcal{K})^3}\nonumber\\
&=\frac{2r^2\mathcal{E}({\mathcal{E}}^2-{r'}^2{\mathcal{K}}^2)-2{r'}^2(\mathcal{K}-\mathcal{E})^2(\mathcal{E}-{r'}^2\mathcal{K})}
{r{r'}^2(\mathcal{E}-{r'}^2\mathcal{K})^3}\nonumber\\
&=\frac{2r\mathcal{E}({\mathcal{E}}^2-{r'}^2{\mathcal{K}}^2)}{{r'}^2(\mathcal{E}-{r'}^2\mathcal{K})^3}\left[1-g_{1}(r)\right]
\end{align}
where
\begin{equation}
g_{1}(r)=\left[\frac{{r'}^{3/4}(\mathcal{K}-\mathcal{E})}{r^2}\right]^2\cdot\frac{{r'}^{1/2}}{\mathcal{E}}\cdot
\frac{r^2}{\mathcal{E}-{r'}^2\mathcal{K}}\cdot
\frac{(\mathcal{E}-{r'}^2\mathcal{K})^2}{{\mathcal{E}}^2-{r'}^2{\mathcal{K}}^2}.
\end{equation}

It follows from Lemma 2.2(1), (2), (5) and (6) together with  (2.4) that $g_{1}(r)$ is strictly decreasing from $(0,1)$ onto $(0,1)$. Then equation (2.3) leads to the conclusion that $g(r)$ is strictly increasing in $(0,1)$. Making use of l'H\^{o}pital's rule one has
\begin{align*}
\lim_{r\rightarrow 0^+}g(r)=&\lim_{r\rightarrow
0^+}\frac{(\mathcal{K}-\mathcal{E})(\mathcal{E}-{r'}^2\mathcal{K})/r^4+\mathcal{E}[(\mathcal{K}-\mathcal{E})-(\mathcal{E}-{r'}^2\mathcal{K})]/r^4}
{(\mathcal{E}-{r'}^2\mathcal{K})^2/r^4}\\
=&\frac{{\pi}^2/16+\pi/2\times \pi/16}{{\pi}^2/16}=\frac{3}{2},
\end{align*}
while $g(1^-)=+\infty$ is clear. $\Box$

\medskip
{\bf Lemma 2.4.} Let $r\in(0,1)$, $p\geq 1/2$ and
\begin{equation*}
h(r)=(2p-1)r^2+2p\frac{r^2\mathcal{E}}{\mathcal{E}-{r'}^2\mathcal{K}}.
\end{equation*}
Then $h(r)$ is strictly decreasing from $(0,1)$ onto $(4p-1,4p)$ if
and only if $p\leq 2$.

\medskip
{\bf Proof.} Differentiating $h(r)$ leads to
\begin{align}
h'(r)=&2(2p-1)r+2p\frac{r(3\mathcal{E}-\mathcal{K})(\mathcal{E}-{r'}^2\mathcal{K})-r^3\mathcal{K}\mathcal{E}}{(\mathcal{E}-{r'}^2\mathcal{K})^2}\nonumber\\
=&2pr\left[2-\frac{1}{p}-g(r)\right],
\end{align}
where $g(r)$ is defined as in Lemma 2.3.

Therefore, Lemma 2.4 follows easily from Lemma 2.3 and (2.5) together with the limiting values $h(0^+)=4p$ and $h(1^{-})=4p-1$. $\Box$

\medskip
{\bf Lemma 2.5.} The inequality
\begin{equation*}
\frac{1}{4p}<\left(\frac{4}{\pi}\right)^{1/p}-1<\frac{1}{4p-1}
\end{equation*}
holds for $p\in[1/2,2]$.

\medskip
{\bf Proof.} Let $f_{1}(p)=\left[(4p+1)/(4p)\right]^{p}$ and
$f_{2}(p)=\left[4p/(4p-1)\right]^p$. Then it is not difficult to
verify that $f_{1}(p)$ is strictly increasing in $[1/2,2]$ and
$f_{2}(p)$ is strictly decreasing in  $[1/2,2]$. Hence, we get
\begin{equation}
f_{1}(p)<f_{1}(2)=\frac{81}{64}=1.265625\cdots<\frac{4}{\pi}=1.273238\cdots
\end{equation}
and
\begin{equation}
f_{2}(p)>f_{2}(2)=\frac{64}{49}=1.306122\cdots>\frac{4}{\pi}=1.273238\cdots
\end{equation}
for $p\in[1/2,2]$.

Therefore, Lemma 2.5 follows from (2.6) and (2.7). $\Box$

\medskip
{\bf Lemma 2.6.} Let $u\in[0,1]$, $p\in[1/2,2]$, $r\in(0,1)$ and
\begin{equation}
f_{u,p}(r)=p\log(1+ur^2)-\log\left[\frac{2}{\pi}(2\mathcal{E}-{r'}^2\mathcal{K})\right].
\end{equation}
Then $f_{u,p}(r)<0$ for $r\in(0,1)$ if and only if $u\leq 1/(4p)$, and
$f_{u,p}(r)>0$ for $r\in(0,1)$ if and only if $u\geq (4/\pi)^{1/p}-1$.

\medskip
{\bf Proof.} From (2.8) we get
\begin{equation}
f_{u,p}(0^+)=0,
\end{equation}
\begin{equation}
f_{u,p}(1^-)=p\log(1+u)+\log\left(\frac{\pi}{4}\right)
\end{equation}
and
\begin{align}
{f_{u,p}}'(r)=&\frac{2pur}{1+ur^2}-\frac{\mathcal{E}-{r'}^2\mathcal{K}}{r(2\mathcal{E}-{r'}^2\mathcal{K})}\nonumber\\
=&\frac{(2p-1)r^2(\mathcal{E}-{r'}^2\mathcal{K})+2pr^2\mathcal{E}}{r(1+ur^2)(2\mathcal{E}-{r'}^2\mathcal{K})}\left[u-\frac{1}{h(r)}\right],
\end{align}
where $h(r)$ is defined as in Lemma 2.4.

It follows from  Lemma 2.4 that $1/h(r)$ is strictly increasing from $(0,1)$ onto $(1/(4p),1/(4p-1))$ if $1/2\leq p\leq 2$.

Making use of Lemma 2.5, we divide the proof into four cases.

{\bf Case 1} $u\leq 1/(4p)$. Then (2.11) and the monotonicity of $1/h(r)$ lead to the conclusion that $f_{u,p}(r)$ is strictly decreasing in $(0,1)$. Therefore, $f_{u,p}(r)<0$ for $r\in(0,1)$ follows from (2.9) and the monotonicity of $f_{u,p}(r)$.

{\bf Case 2} $u\geq 1/(4p-1)$. Then (2.11) and the monotonicity of $1/h(r)$ lead to the conclusion that $f_{u,p}(r)$ is strictly increasing in $(0,1)$. Therefore, $f_{u,p}(r)>0$ for $r\in(0,1)$ follows from (2.9) and the monotonicity of $f_{u,p}(r)$.

{\bf Case 3} $(4/\pi)^{1/p}-1<u<1/(4p-1)$. Then from (2.10) and (2.11) together with the monotonicity of $1/h(r)$ we clearly see that there exists $\lambda\in(0,1)$ such that $f_{u,p}(r)$ is strictly increasing in $(0,\lambda]$ and strictly decreasing in $[\lambda,1)$, and
\begin{equation}
f_{u,p}(1^-)\geq 0.
\end{equation}

Therefore, $f_{u,p}(r)>0$ for $r\in(0,1)$ follows from (2.9) and (2.12) together with the piecewise monotonicity of $f_{u,p}(r)$.

{\bf Case 4} $1/(4p)<u<(4/\pi)^{1/p}-1$. Then from (2.10) and (2.11) together with the monotonicity of $1/h(r)$ we clearly see that there exists $0<\mu<1$ such that $f_{u,p}(r)$ is strictly increasing in $(0,\mu]$ and strictly decreasing in $[\mu,1)$, and
\begin{equation}
f_{u,p}(1^-)<0.
\end{equation}

Therefore, equation (2.9) and inequality (2.13) together with the piecewise monotonicity of $f_{u,p}(r)$ lead to the conclusion that there exists $0<\mu<\eta<1$ such that $f_{u,p}(r)>0$ for $r\in(0,\eta)$ and $f_{u,p}(r)<0$ for $r\in(\eta,1)$. $\Box$

\medskip
{\bf Lemma 2.7.} If $r\in(0,1)$, then the function
\begin{equation*}
F(r)=(2{\mathcal{E}}-{r'}^2{\mathcal{K}})^2\left[1+\frac{{\pi}^2-4(2{\mathcal{E}}-{r'}^2{\mathcal{K}})^2}{{\pi}^2 r^2}\right]
\end{equation*}
is strictly increasing from $(0,1)$  onto
$({\pi}^2/8,8({\pi}^2-8)/{\pi}^2)$.

\medskip
{\bf Proof.} Differentiating $F(r)$ leads to
\begin{align}
f'(r)=&\frac{2(2{\mathcal{E}}-{r'}^2{\mathcal{K}})({\mathcal{E}}-{r'}^2{\mathcal{K}})}{r}\left[1+\frac{{\pi}^2-4(2{\mathcal{E}}-{r'}^2{\mathcal{K}})^2}{{\pi}^2
r^2}\right]\nonumber\\
&+(2{\mathcal{E}}-{r'}^2{\mathcal{K}})^2\frac{8\mathcal{E}(2{\mathcal{E}}-{r'}^2{\mathcal{K}})-2{\pi}^2}{{\pi}^2
r^3}\nonumber\\
=&\frac{2{r'}^2(2{\mathcal{E}}-{r'}^2{\mathcal{K}})({\mathcal{E}}-{r'}^2{\mathcal{K}})}{r^3}\left[\frac{4(2{\mathcal{E}}-{r'}^2{\mathcal{K}})^2-{\pi}^2}{{\pi}^2r^2}\frac{r^2{\mathcal{K}}}{({\mathcal{E}}-{r'}^2{\mathcal{K}})}-1\right].
\end{align}

It follows from  Lemma 2.2(4) and (7) together with (2.14) that $F(r)$ is
strictly increasing in $(0,1)$. Moreover, the limiting values
$F(0^+)={\pi}^2/8$ and $F(1^-)=8({\pi}^2-8)/{\pi}^2$ follows easily from Lemma 2.2(7). $\Box$

\bigskip
\bigskip
\noindent{\bf\large 3. Proof of Theorems 1.1 and 1.2}
\bigskip
\setcounter{section}{3} \setcounter{equation}{0}

{\bf Proof of Theorem 1.1.} Let $q\in(0,1/2]$ and
$t=(1-r')/(1+r')\in(0,1)$. Then
\begin{equation}
{\mathcal{E}}(r)={\mathcal{E}}\left(\frac{2\sqrt{t}}{1+t}\right)=\frac{2{\mathcal{E}}(t)-{t'}^2{\mathcal{K}}(t)}{1+t}
\end{equation}
and
\begin{align}
&\sqrt{q+(1-q){r'}^2}+\sqrt{(1-q)+q{r'}^2}\nonumber\\
=&\frac{\sqrt{1+t^2+2(1-2q)t}+\sqrt{1+t^2-2(1-2q)t}}{1+t}.
\end{align}

From Lemma 2.2(7) we get
\begin{equation}
\frac{8}{{\pi}^2}\left[2{\mathcal{E}}(t)-{t'}^2{\mathcal{K}}(t)\right]^2>1+t^2.
\end{equation}

Equations (3.1) and (3.2) together with inequality (3.3) lead to the following equivalence relations:
\begin{align}
&\frac{4}{\pi}{\mathcal{E}}(r)>\sqrt{q+(1-q){r'}^2}+\sqrt{(1-q)+q{r'}^2}\nonumber\\
\Longleftrightarrow&\frac{4}{\pi}\left[2{\mathcal{E}}(t)-{t'}^2{\mathcal{K}}(t)\right]>\sqrt{1+t^2+2(1-2q)t}+\sqrt{1+t^2-2(1-2q)t}\nonumber\\
\Longleftrightarrow&\frac{8}{\pi^2}\left[2{\mathcal{E}}(t)-{t'}^2{\mathcal{K}}(t)\right]^2-(1+t^2)>\sqrt{(1+t^2)^2-4(1-2q)^2t^2}\nonumber\\
\Longleftrightarrow&\frac{64}{\pi^4}\left[2{\mathcal{E}}(t)-{t'}^2{\mathcal{K}}(t)\right]^4-\frac{16}{{\pi}^2}(1+t^2)\left[2{\mathcal{E}}(t)-{t'}^2{\mathcal{K}}(t)\right]^2
>-4(1-2q)^2t^2\nonumber\\
\Longleftrightarrow&\frac{4}{\pi^2}\left[2{\mathcal{E}}(t)-{t'}^2{\mathcal{K}}(t)\right]^2\left[1+\frac{{\pi}^2-4(2{\mathcal{E}}(t)-{t'}^2{\mathcal{K}}(t))^2}{{\pi}^2t^2}\right]
<(1-2q)^2.
\end{align}

Similarly, we have
\begin{align}
&\frac{4}{\pi}{\mathcal{E}}(r)<\sqrt{q+(1-q){r'}^2}+\sqrt{(1-q)+q{r'}^2}\Longleftrightarrow\nonumber\\
&\frac{4}{\pi^2}\left[2{\mathcal{E}}(t)-{t'}^2{\mathcal{K}}(t)\right]^2
\left[1+\frac{{\pi}^2-4(2{\mathcal{E}}(t)-{t'}^2{\mathcal{K}}(t))^2}{{\pi}^2t^2}\right]>(1-2q)^2.
\end{align}

Therefore, Theorem 1.1  follows easily from (3.4) and (3.5) together with Lemma 2.7. $\Box$

\bigskip
In order to prove Theorem 1.2, we need to introduce the Toader mean
which is closely related to the complete elliptic integral of the
second kind. For $a,b>0$ with $a\neq b$, the Toader mean $T(a,b)$
[17, 18, 20, 25] is defined by
\begin{align}
T(a,b)&=
\frac{2}{\pi}\int_{0}^{{\pi}/{2}}\sqrt{a^2{\cos^2{\theta}}+b^2{\sin^2{\theta}}}d\theta=\left\{\begin{array}{ll}
{2a}{\mathcal{E}}\big{(}{\sqrt{1-(b/a)^2}}\big{)}/{\pi},&a>b,\\
{2b}{\mathcal{E}}\big{(}{\sqrt{1-(a/b)^2}}\big{)}/{\pi},&a<b.
\end{array}\right.
\end{align}
In particular,
\begin{equation}
\mathcal{E}(r)=\pi T(1,r')/2
\end{equation}
for $r\in(0,1)$.

\medskip
{\bf Proof of Theorem 1.2.} Let $a,b>0$ with $a\neq b$, $t, t_{1}, t_{2}\in[1/2,1]$, $p\in[1/2,2]$ and
\begin{equation}
Q_{t,p}(a,b)=C^{p}(ta+(1-t)b, tb+(1-t)a)A^{1-p}(a,b),
\end{equation}
where $A(a,b)=(a+b)/2$ and $C(a,b)=(a^2+b^2)/(a+b)$ are the arithmetic and contraharmonic means of $a$ and $b$, respectively. Then $Q_{t,p}(a,b)$ is strictly increasing with respect to $t\in[1/2,1]$ for fixed $a,b>0$ with $a\neq b$.

Since both $Q_{t,p}(a,b)$ and $T(a,b)$ are symmetric and homogeneous
of degree 1. Without loss of generality, we assume that $a>b$. Let
$x=(a-b)/(a+b)\in(0,1)$. Then (3.6)-(3.8) lead to
\begin{align}
\log\left(\frac{Q_{t,p}(a,b)}{T(a,b)}\right)=&\log\left(\frac{Q_{t,p}(a,b)}{A(a,b)}\right)-\log\left(\frac{T(a,b)}{A(a,b)}\right)\nonumber\\
=&p\log\left[1+(1-2t)^2x^2\right]-\log\left[\frac{T(1,b/a)}{A(1,b/a)}\right]\nonumber\\
=&p\log\left[1+(1-2t)^2x^2\right]-\log\left[\frac{2}{\pi}(1+x)\mathcal{E}\left(\frac{2\sqrt{x}}{1+x}\right)\right]\nonumber\\
=&p\log\left[1+(1-2t)^2x^2\right]-\log\frac{2}{\pi}\left[2\mathcal{E}(x)-{x'}^2\mathcal{K}(x)\right].
\end{align}

From Lemma 2.6 and equation (3.9) we clearly see that the double inequality
\begin{equation}
Q_{t_{1},p}(a,b)<T(a,b)<Q_{t_{2},p}(a,b)
\end{equation}
holds if and only if $t_{1}\leq 1/2+\sqrt{1/(4p)}/2 $ and $t_{2}\geq
1/2+\sqrt{(4/\pi)^{1/p}-1}/2$.

In particular, if $a=1$ and $b=r'$, then inequality (3.10) becomes
\begin{align}
&2^{p-2}\pi(1+r')^{1-2p}\left\{\left[t_{1}+(1-t_{1})r'\right]^2+\left[(1-t_{1})+t_{1}r'\right]^2\right\}^{p}<\mathcal{E}(r)\nonumber\\
&<2^{p-2}\pi(1+r')^{1-2p}\left\{\left[t_{2}+(1-t_{2})r'\right]^2+\left[(1-t_{2})+t_{2}r'\right]^2\right\}^{p}.\nonumber
\end{align}

The proof of Theorem 1.2 is therefore complete. $\Box$

\medskip
Let $p=2$ and $p=1/2$ in the first and second inequality of (1.5),
respectively. Then we get

\medskip
{\bf Corollary 3.1.} The double inequality
\begin{align*}
\frac{\pi}{(1+r')^3}\left\{\left[\lambda+(1-\lambda)r'\right]^2+\left[(1-\lambda)+\lambda r'\right]^2\right\}^{2}
<\mathcal{E}(r)\\
<\frac{\pi}{2}\left\{\frac{\left[\mu+(1-\mu)r'\right]^2+\left[(1-\mu)+\mu r'\right]^2}{2}\right\}^{1/2}
\end{align*}
holds for all $r\in(0,1)$ with $\lambda=1/2+\sqrt{2}/8$ and $\mu=1/2+\sqrt{(4/\pi)^2-1}/2$.

\bigskip
\bigskip
\noindent{\bf\large 4. Comparison of the bounds for $\mathcal{E}(r)$}
\bigskip
\setcounter{section}{4} \setcounter{equation}{0}

In this section, we compare our bounds for $\mathcal{E}(r)$ with the
bounds in (1.1)-(1.3).

\medskip
{\bf Remark 4.1.} Simple computations lead to that the upper bound
in (1.4) is exactly the bound in (1.3).

\medskip
{\bf Remark 4.2.} Let $x\in(0,1)$ and $\mu=1/2+\sqrt{(4/\pi)^2-1}/2$. Then simple computation leads to
\begin{equation*}
(1+x^2)-\left\{[\mu+(1-\mu)x]^2+[(1-\mu)+\mu
x]^2\right\}=(1-\frac{8}{{\pi}^2})(1-x)^2>0.
\end{equation*}

Therefore, the upper bound in Corollary 3.1 is better than the bound
in (1.2).

\medskip
{\bf Remark 4.3.} Let $\overline{\alpha}=1/2-\sqrt{2}/4$, $\overline{\beta}=1/2+\sqrt{2}/4$ and $\mu=1/2+\sqrt{(4/\pi)^2-1}/2$.
Then we have
\begin{align}
&\lim\limits_{r\rightarrow
1}\left[\frac{\pi}{4}\left(\sqrt{1-\overline{\alpha}r^2}+\sqrt{1-\overline{\beta}r^2}\right)\right]\nonumber\\
=&\frac{\pi}{8}\left(\sqrt{2+\sqrt{2}}+\sqrt{2-\sqrt{2}}\right)=1.026172\cdots
\end{align}
and
\begin{equation}
\lim\limits_{r\rightarrow
1}\frac{\pi}{2}\left\{\frac{\left[\mu+(1-\mu)r'\right]^2+\left[(1-\mu)+\mu r'\right]^2}{2}\right\}^{1/2}=1.
\end{equation}

Equations (4.1) and (4.2) show that the upper bound in Corollary 3.1
is asymptotically precise when $r\rightarrow 1$, and there exists
$0<\delta_{1}<1$ such that the upper bound in Corollary 3.1 is
better than the bound in (1.3) when $r\in(1-\delta_{1},1)$.

\medskip
{\bf Remark 4.4.} Let $\beta=1/2-2\sqrt{2({\pi}^2-8)}/{\pi}^2$, then
\begin{align}
&\lim\limits_{r\rightarrow
1}\frac{\pi}{4}\left(\sqrt{\beta+(1-\beta){r'}^2}+\sqrt{(1-\beta)+\beta{r'}^2}\right)=1
\end{align}
and
\begin{equation}
\lim\limits_{r\rightarrow
1}\frac{\pi}{2}\left(\frac{1+{r'}^{3/2}}{2}\right)^{2/3}=2^{-5/3}\pi=0.989539\cdots.
\end{equation}

Equations (4.3) and (4.4) imply that the lower bound in Theorem 1.1
with $\beta=1/2-2\sqrt{2({\pi}^2-8)}/{\pi}^2$ is asymptotically
precise when $r\rightarrow 1$, and there exists $0<\delta_{2}<1$
such that the lower bound in Theorem 1.1 is better than the bound in
(1.1) when $r\in(1-\delta_{2},1)$.

\medskip
{\bf Remark 4.5.} Let $\lambda=1/2+\sqrt{2}/8$, $r\in(0,1)$ and
$x=(1-r^2)^{1/4}\in(0,1)$. Then the following equivalence relations
lead to the conclusion that the lower bound in Corollary 3.1 is
better than the bound in inequality (1.1).
\begin{align}
&\frac{\pi}{(1+r')^3}\left\{\left[\lambda+(1-\lambda)r'\right]^2+\left[(1-\lambda)+\lambda r'\right]^2\right\}^{2}>\frac{\pi}{2}\left(\frac{1+{r'}^{3/2}}{2}\right)^{2/3}\nonumber\\
\Longleftrightarrow&\frac{\left[9+9(1-r^2)+14\sqrt{1-r^2}
\right]^2}{128(1+\sqrt{1-r^2})^3}>
\left[\frac{1+(1-r^2)^{3/4}}{2}\right]^{2/3}\nonumber\\
\Longleftrightarrow&\frac{\left(9+9x^4+14x^2\right)^2}{128(1+x^2)^3}>
\left(\frac{1+x^3}{2}\right)^{2/3}\nonumber\\
\Longleftrightarrow&(9+9x^4+14x^2)^6-524288(1+x^2)^9(1+x^3)^2> 0\nonumber\\
\Longleftrightarrow&(x-1)^4\big[7153x^{20}+28612x^{19}+313054x^{18}+60580x^{17}+2074909x^{16}\nonumber\\
&-277424x^{15}+6613736x^{14}-1390192x^{13}+12597746x^{12}-2615880x^{11}\nonumber\\
&+15507060x^{10}-2615880x^{9}+12597746x^{8}-1390192x^{7}+6613736x^{6}\nonumber\\
&-277424x^{5}+2074909x^{4}+60580x^{3}+313054x^{2}+28612x+7153\big]>0.\nonumber
\end{align}

\bigskip
\noindent{\bf\normalsize Acknowledgements}
\medskip

This research was supported by the Natural Science Foundation of
China under Grant 11071069, and the Innovation Team Foundation of
the Department of Education of Zhejiang Province under Grant
T200924.

\bigskip
\def\refname
\end{document}